\documentclass[11pt]{article}

\usepackage{latexsym,amssymb, amsmath}
\usepackage{comment,theorem}
\usepackage[colorlinks=true,hyperindex=true]{hyperref}

\makeatletter \@addtoreset{equation}{section}

\newtheorem{thm}{Theorem}[section]
\newtheorem{prop}[thm]{Proposition}
\newtheorem{lem}[thm]{Lemma}
\newtheorem{cor}[thm]{Corollary}
\newtheorem{conj}[thm]{Conjecture}

\theorembodyfont{\rmfamily}

\newtheorem{rema}[thm]{Remark}

\def\qed{\hfill{$\blacksquare$}\medskip}

\def\bp{\noindent{\it Proof.}\ }

\def\2#1{{\cal #1}}
\def\6#1{{\mathfrak #1}}
\def\7#1{{\mathbb #1}}

\newcommand{\End}{\mathrm{End}}

\newcommand{\Tr}{{\rm Tr}}

\newcommand{\rarr}{\rightarrow}

\newcommand{\nn}{\nonumber}

\newcommand{\be}{\begin{equation}}
\newcommand{\ee}{\end{equation}}
\newcommand{\bean}{\begin{eqnarray*}}
\newcommand{\eean}{\end{eqnarray*}}
\newcommand{\bea}{\begin{eqnarray}}
\newcommand{\eea}{\end{eqnarray}}

\begin{document}

\title{An integral formula for Lie groups, and the Mathieu conjecture reduced to Abelian non-Lie conjectures}

\author{M.\ M\"uger\footnote{Radboud University, Nijmegen, The Netherlands. {\tt mueger@math.ru.nl}} \ and L.\
Tuset\footnote{Oslo Metropolitan University, Oslo, Norway. \tt{larst@oslomet.no}}}
\date{\today}
\maketitle

\begin{abstract}
We present an explicit integration formula for the Haar integral on a compact connected Lie group. This
formula relies on a known decomposition of a compact connected simple Lie group into symplectic leaves, when
one views the group as a Poisson Lie group. In this setting the Haar integral is constructed using the
Kostant harmonic volume form on the corresponding flag manifold, and explicit coordinates are known for these 
invariant differential forms. The formula obtained is related to one found by Reshetikhin-Yakimov.  

Using our integration formula, we reduce the Mathieu conjecture to two stronger conjectures about Laurent
polynomials in several complex variables with polynomial coefficients in several real variable polynomials. In
these stronger conjectures there is no reference to Lie group theory. 
\end{abstract}

\section{Introduction}

Recall that the Mathieu conjecture \cite{mathieu} says that for any compact connected Lie group $K$, and for
any pair $f,g$ of $K$-finite complex valued functions, it is true that if $\int_K f^n (k)\, dk =0$ for all
$n\in\mathbb{N}$ with respect to the Haar measure $dk$, then $\int_K (f^n g)(k)\, dk =0$ for all but finitely
many $n\in\mathbb{N}$.  

This conjecture implies \cite{mathieu} the famous Jacobian conjecture. A difficult proof involving algebraic
geometry \cite{duist-kallen} shows that the conjecture holds when $K$ is a compact torus, but beyond this, nothing
substantial seems to be known.  The point of the present paper is to reduce the Mathieu conjecture for general
$K$ to a conjecture of a more Abelian nature. This conjecture involves complex valued functions with no
reference to Lie groups and their complicated structure and representation theory. It can be considered as a
generalization of the abelian case of the Mathieu conjecture, known to be true, in a Lie-free direction.
We also state a second, more natural, conjecture implying our first one, which again generalizes a statement 
proven in \cite{duist-kallen}. 

All this hinges on an explicit formula for the Haar integral on $K$. Let us briefly explain how this formula
is obtained. By a known structure theorem, $K$ is the product of a compact torus and finitely many compact
connected simple Lie groups, modulo a finite subgroup of the center of the product. The Haar integral of $K$
correspondingly reduces essentially to those of the factors, so that we are left with studying the Haar
integral of a compact connected simple Lie group $K$. The latter is a compact real form of a complex simple
Lie group $G$, the complexification $K_\7C$ of $K$. Picking a maximal torus $T$ for $K$, we then consider the
flag manifold $K/T$, and show that if we have a $K$-invariant volume form there, we can use that to produce
the Haar integral on $K$. This is fairly standard. There is a beautiful $K$-invariant volume form on $K/T$,
namely the so called Kostant harmonic form associated to the biggest Schubert cell $\Sigma_{w_0}$; these cells
form a $CW$-complex for $K/T$ which corresponds to the Bruhat decomposition of $G$. The real dimension of
$\Sigma_{w_0}$ is twice the length of the longest element $w_0$ in the Weyl group of $K$. Let $G=KAN$ be the
Iwasawa decomposition of $G$, allowing $G$ to act on $K$ by letting $g\circ k$ be 
the projection of $gk$ onto $K$. Letting $\dot{w}_0\in K$ be a representative of $w_0$, one can show that the
orbit $N\circ\dot{w}_0\subset K$ is diffeomorphic to $\Sigma_{w_0}$. In fact, this orbit is the orbit through
$\dot{w}_0$ of the so called dressing action associated to $K$ seen as a Poisson Lie group equipped with the
standard Poisson bracket. These orbits, one for each representative of each element in the Weyl
group, are the symplectic leaves of the Poisson Lie group. Using the reduced expression of $w_0$ into
fundamental reflections corresponding to the simple roots of $\mathfrak g$, together with the group
homomorphisms $SU(2)\to K$ corresponding to the nodes in the Dynkin diagram of $\mathfrak g$, one obtains a
further decomposition of $N$ and of the symplectic leaves into the ones occurring for $SU(2)$, see also
\cite{soib-vaks}. This yields a very explicit formula for the Kostant harmonic volume form. In fact, we get
coordinate descriptions of various useful quantities, not only of all the Kostant harmonic differential forms,
but for the symplectic differential $2$-forms on each leaf, and their Liouville volume forms, as well as for
the Haar measure $dn$ on $N$, and for some interesting functions $a_w$ that also appear in the work of
Reshetikhin-Yakimov \cite{resh-yak}. In that paper, a formula for the Haar integral on a compact connected
simple Lie group is given in terms of the 
functions $a_{w_0}$ and the Liouville volume form on the biggest leaf. Of course, also the decomposition of
this big cell into $2$-dimensional leaves is given, but no proof is provided for their integration formula,
nor any explicit coordinates, which we really need. We are pretty sure though, that Reshetikhin-Yakimov were
aware of these coordinates, which were exhibited in the crucial article \cite{lu}, and they even refer to this
article in their work. However, their focus was different, namely on quantum groups, and the relation to the
celebrated work in \cite{soib-vaks}.  

We provide all the necessary results from this theory in Sections \ref{sec2} and \ref{s-Poisson}, concentrating on
the parts that are particularly relevant for us. 

In Section \ref{sec-int}, working in the coordinates described in \cite{lu}, we arrive quite rapidly at our 
explicit integration formula for the Haar integral on a compact connected Lie group $K$. We explain carefully
how the formula works, rendering it hopefully useful for people that need some very concrete stuff in
calculations involving the Haar integral. These results hold for all measurable functions.

In Section \ref{sec5} we specialize to the integration formula to finite type functions, setting the stage for
the final Section \ref{sec6}, where we turn to the Mathieu conjecture, rephrasing it using our integration
formula. Here we formulate our own conjectures. The first one is very similar to the Mathieu conjecture,
except that we have dispensed with everything relating to Lie theory. The second one is an even stronger
conjecture, and we explain with an easy proposition why the latter conjecture is indeed stronger. In passing,
we discuss the $SU(2)$-case that was already considered in \cite{mueger-tus}.


\section{Some Lie theory}\label{sec2}
For later reference, we collect a number of results, most of them well known, on Lie groups.

Recall that a compact torus by definition is a compact connected Abelian Lie group. The following result is well known:

\begin{thm} \label{thm-prod}
Every compact connected Lie group $K$ admits an isomorphism $K\cong(K_1\times\cdots\times K_n\times T)/D$,
where the $K_i$ are compact connected simple Lie groups, $T$ is a compact torus (of dimension $\dim Z(K)$),
and $D$ is a finite subgroup of the center of the product.
\end{thm}

For example, $U(2)$ is isomorphic to the quotient of $SU(2)\times\7T^1$ by the normal subgroup
$\{({\bf 1},1),(-{\bf 1},-1)\}$.
For proofs, see e.g.\  \cite[Theorem 6.4.2]{price} or \cite[Section 10.7.2, Theorem 4]{procesi}. 

\begin{thm}\label{th-complexif}
Every compact connected Lie group $K$ admits a complexification $K_\7C$ or $G$, i.e.\ a connected complex Lie group such that:
\begin{itemize}
\item The Lie algebra $\6g$ is isomorphic to the complexification $\6k_\7C$ of $\6k$. (I.e.\ $K$ is a real form of $G$.) 
\item The inclusion map $K\rarr G$ induced by $\6k\rarr\6k_\7C\rarr\6g$ induces an isomorphism $\pi_1(K)\rarr\pi_1(G)$.
\end{itemize}
\end{thm}

For a good presentation see \cite[Chapter 24]{bump}. The same reference gives good discussions of the Iwasawa
(=KAN) and Bruhat decomposition for the reductive group $G$ thus obtained, cf.\ Chapters 26 and 27,
respectively. 

Finally we need a result on integration over a compact connected Lie group, for which we include a proof for
lack of a convenient reference. 

When $X$ is a compact Hausdorff space, let $C(X)$ denote the unital C$^*$-algebra of continuous complex valued
functions on $X$ under pointwise $*$-algebra operations and uniform norm. When $X$ is furthermore a (real)
analytic manifold, denote by $A(X)\subset C^\infty (X)$ the unital $*$-subalgebras of $C(X)$ consisting of
analytic and smooth functions, respectively. By the Stone-Weierstrass theorem, see Theorem 4.3.4. in
\cite{gert}, we know that $A(X)$ is dense in $C(X)$ since already polynomials separate points in $X$.  

If $X$ is a compact connected Lie group, the complex linear span $\mathbb{C}[X]$ of matrix coefficients on $X$
is a unital $*$-subalgebra of $A(X)$ which is dense in $C(X)$. Recall that a matrix coefficient on $X$ is a
function of the form $x\mapsto (U(x)u|v)$, where $U\colon X\to B(V)$ is a unitary representation of $X$ on a
finite dimensional Hilbert space $V$ with inner product $(\cdot |\cdot )$, which we take to be conjugate
linear in the second variable. Since $U$ is continuous and $X$ is compact, we know that $U(X)$ is compact, and
thus a closed subgroup of $GL (V)$, so it is a Lie group by Theorem 2.12.6 in \cite{var}. Furthermore, $U$ is
automatically an analytic homomorphism between the Lie groups $X$ and $U(X)$ due to Theorem 2.6 in
\cite{helgason}. Hence $\mathbb{C}[X]\subset A(X)$, and denseness of $\mathbb{C}[X]$ in $C(X)$ is clear from
the Stone-Weierstrass theorem and Corollary 4.22 in \cite{knapp}. One refers to elements of $\mathbb{C}[X]$ as
$X$-finite functions. 

\begin{lem} \label{lem-int}
Let $K$ be a compact topological group with a closed subgroup $T$, and let 
$$
C^T (K)=\{f\in C(K)\, |\, f(kt)=f(k)\, \, \, \forall k\in K, t\in T\}.
$$   
Consider the compact Hausdorff space $K/T$ of left cosets $kT$, and let $\pi\colon K\to K/T$ be the quotient
map $\pi (k)=kT$, which is both continuous and open. The map $\Phi\colon C(K/T)\to C^T (K)$ given by
$\Phi(h)=h\circ\pi$ is a unital $*$-isomorphism between C$^*$-algebras.  

Let $dt, dk$ be the normalized Haar measures on $T$ and $K$, respectively, and let $d(kT)$ be a normalized
$K$-invariant regular Borel measure on $K/T$, which is a homogeneous space under the action of $ K$ by left
translation. The formula $P(f)(k)=\int_T f(kt)\, dt$ for $f\in C(K)$, defines a surjective positive unital
linear map $P\colon C(K)\to C^T (K)\subset C(K)$ such that $P^2 =P$. Moreover, we have 
$$
\int_K f(k)\, dk =\int_{K/T}\Phi^{-1}P(f)(kT)\, d(kT)=\int_{K/T}\int_T f(kt)\, dt\, d(kT)
$$  
for $f\in C(K)$.

If in addition $K$ is a Lie group, then $T$ is a Lie subgroup of $K$ by Theorem 2.12.6 in \cite{var}, and
$K/T$ is by Theorem 2.9.4 in \cite{var} an analytic manifold with analytic $\pi$ and an analytic action
$K\times K/T\to K/T$ of $K$. Furthermore, this is the unique analytic structure such that $h\in C(K/T)$ is
analytic whenever $h\circ\pi$ is analytic, so $\Phi (A(K/T))$ is the unital $*$-subalgebra of $C^T(K)$
consisting of all $f\in A(K)$ such that $f(kt)=f(k)$ for $k\in K$ and $t\in T$. We have
$P(\mathbb{C}[K])\subset \Phi (A(K/T))$. 

Let $\omega$ be a normalized analytic volume form on $K/T$. Then there exists a normalized (positive) regular
Borel measure $d(kT)$ on $K/T$ such that $\int_{K/T}\ h\omega =\int_{K/T} h(kT)\, d(kT)$ for $h\in C^\infty
(K/T)$. This Borel measure is $K$-invariant if and only if $\omega$ is $K$-invariant, meaning that $L_k^*
(\omega )=\omega$ for all $k\in K$, where $L_k^*$ is the map on differential forms and functions induced from
left translation $L_k$ by $k\in K$ on $K/T$.  
Moreover, if $E$ is a finite union of submanifolds of $K/T$ all of dimension strictly less than
$\mathrm{dim}(K/T)$, then $\int_E h(kT)\, d(kT)=0$ for any continuous function $h$ on $K/T$.   
\end{lem}

\bp
Now $\pi$ is continuous by definition of the quotient topology, and the image under $\pi$ of an open set is open by paragraph 1.4.11 in \cite{gert} as $\bigcup_{t\in T}At$ is open in $K$ for any open $A\subset K$. 
 
We have $\Phi (h)\in C^T (K)$ since $h\circ\pi$ is continuous and $h\circ\pi (kt)=h\circ \pi (k)$ for $k\in K$ and $t\in T$. Clearly $\Phi$ is a unital $*$-homomorphism. With respect to uniform norms we have
$\|h\|_{K/T} =\sup_{kT\in K/T} |h(kT)|=\sup_{k\in K}|h\circ\pi (k)|=\|h\circ \pi\|_K$, so $\Phi$ is isometric, and thus injective. It remains to show that it is surjective. Let $f\in C^T (F)$. Define $h\colon K/T\to\mathbb{C}$ by $h\circ\pi =f$. This is well-defined because if $\pi (k')=\pi (k)$, then $k'=kt$ for some $t\in T$, and therefore $h(\pi (k'))=f(k')=f(kt)=f(k)=h(\pi (k))$. By definition of the topology on $K/T$ as the final topology induced by $\pi$, we know that $h$ is continuous as $f$ is, see Proposition 1.4.10 in \cite{gert}, 
so $h\in C(K/T)$ and $\Phi (h)=f$. We have shown that $\Phi$ is a $*$-isomorphism.
 
Alternatively, to see that the last $h$ above is continuous, consider open $U\subset\mathbb{C}$. Then 
$$
h^{-1}(U)=\{ kT\in K/T\, |\, h(kT)\in U\, \, \forall k\in K\}
=\{\pi (k)\, |\, (h\circ\pi )(k)\in U\, \, \forall k\in K\}
$$
and $(h\circ\pi )(k)\in U$ means that $k\in (h\circ\pi )^{-1} (U)=f^{-1}(U)$, so 
$h^{-1}(U)=\pi (f^{-1}(U))$. Now $f^{-1}(U)$ is open as $f$ is continuous, and $\pi$ is open, so $h^{-1}(U)$ is open, and $h$ is continuous. 

Next, let $f\in C(K)$. We claim that $P(f)\in C^T (K)$. Recall that a continuous function on a compact group
is uniformly continuous in the sense that for every $\varepsilon>0$ there is an open neighborhood $U$ of the
identity such that $|f(g)-f(h)|<\varepsilon$ whenever $gh^{-1}\in U$. Since $(gt)(ht)^{-1}=gh^{-1}\in U$, we
have 
\[ |P(f)(g)-P(f)(h)|=\Big|\int_T (f(gt)-f(ht))dt \Big|\le\int_T\varepsilon\,dt=\varepsilon.\]
Furthermore, $P(f)$ belongs to $C^T (K)$ because by $T$-invariance of $dt$, we get
$kt'\mapsto \int_T f(kt't)\, dt =\int_T f(kt)\, dt$ for $t'\in T$. If already $f\in C^T (K)$, then
$P(f) (k)=\int_T f(kt)\, dt=\int_T f(k)\, dt=f(k)$, so $P(f)=f$, showing that $P$ is surjective and that it is
idempotent, i.e. satisfies $P^2 =P$. Clearly $P$ is unital and positive.  

This leads to the next identity in the lemma because we have a unital positive linear functional on $C(K)$ given by $f\mapsto\int_{K/T}\Phi^{-1}P(f)(kT)\, d(kT)$ such that 
$$
\int_{K/T}\Phi^{-1}PL^*_{k'}(f)(kT)\, d(kT)=
\int_{K/T}\Phi^{-1}P(f)(k'kT)\, d(kT)=\int_{K/T}\Phi^{-1}P(f)(kT)\, d(kT)
$$
for $k'\in K$ by $K$-invariance of $d(kT)$. So it is therefore the normalized Haar integral on $K$ by uniqueness of such integrals, see Theorem 6.6.12 in \cite{gert}. Letting $h=\Phi^{-1}P(f)$, we get $P(f)=\Phi (h)=h\circ\pi$, and $h(kT)=h\circ\pi (k)=P(f)(k)=\int_T f(kt)\, dt$, which gives the second integration formula for the Haar integral on $K$.   

Moving on to the Lie group case, to establish the inclusion $P(\mathbb{C}[K])\subset \Phi (A(K/T))$, it suffices to show that $P(f)\in A(K)$ for any $f\in \mathbb{C}[K]$, and by linearity we may assume that $f$ is a matrix coefficient, say $f(k)=(U(k)u|v)$ for $k\in K$ and $u,v\in V$. The map $(x,y)\mapsto\int_T (U(t)x|y)\, dt$ is a (bounded) sesquilinear form on $V$, so by Lemma 3.2.2 in \cite{gert} there exists $p\in B(V)$ such that 
$(px|y)=\int_T (U(t)x|y)\, dt$ for all $x,y\in V$. 
Then 
$$
P(f)(k)=\int_T (U(kt)u|v)\, dt =\int_T (U(t)u|U(k^{-1})v)\, dt =(pu|U(k^{-1})v)=(U(k)(pu)|v),
$$  
which shows that $P(f)$ is analytic, being the composition of two analytic maps. 

For the next assertion in the lemma, define a unital positive linear functional $\varphi$ on $C(K/T)$ as
follows. Let $h\in C(K/T)$. Pick a sequence of functions $h_n\in A(K/T)$ such that $h_n\to h$ in the uniform
norm. This can be done by denseness, see the paragraph before this lemma. Set
$\varphi(h)=\lim_{n\to\infty}\int_{K/T} h_n\omega$. We claim that $\varphi$ is well-defined. If we have
another sequence of analytic functions $h'_n$ that also converges to $h$ uniformly, then for any  
$\varepsilon >0$ there exists $N$ such that both $\|h-h_n\|$ and $\|h-h'_n\|$ are less than $\varepsilon$ for 
all $n>N$. Then $\|h_n-h'_n\|<2\varepsilon$ for all $n>N$ by the triangle inequality, so 
$|\int_{K/T} h_n\omega -\int_{K/T} h'_n\omega |\leq\int_{K/T} \| h_n -h'_n\|\omega <2\varepsilon$ for all 
$n>N$, thus yielding the same two numbers $\varphi (h)$. So we get a map $\varphi$, and its claimed properties
clearly hold. By Theorem 2.14 and Theorem 2.17 in \cite{rudin}, we get the required normalized regular Borel
measure $d(kT)$, and the associated claimed identity in the lemma holds by definition of $\varphi$.   

If $\omega$ is $K$-invariant, then for $h\in A(K/T)$ we know that $L^*_k (h)\in A(K/T)$, and by equation 1.1.23 in \cite{var}, we get
$$
\int_{K/T} L^*_k (h) (kT)\, d(kT)=\int_{K/T} L^*_k (h)\omega =\int_{K/T} L^*_k (h)L^*_k (\omega )=\int_{K/T} L^*_k (h\omega ),
$$
which equals
$\int_{K/T} h\omega =\int_{K/T} h(kT)\, d(kT)$ since $L_k$ is an orientation preserving analytic diffeomorphism, see equation 1.1.18 in \cite{var}. The extension to continuous $h$ follows by denseness and Lebesgue's dominated convergence theorem, see Theorem 1.3.4 in \cite{rudin}, together with the fact that $L^*_k$ is continuous with respect to uniform norm. Conversely, if $d(kT)$ is $K$-invariant, then as above we get $\int_{K/T} L^*_k (h)\omega=\int_{K/T} h\omega=\int_{K/T} L^*_k (h)L^*_k (\omega)$ for any smooth function $h$, and by the richness of such functions, see Lemma 1.2 in \cite{helgason}, we get
$L^*_k (\omega) =\omega$.

For the last claim in the lemma, by compactness of $K/T$, we can cover the oriented manifold by finitely many orientation preserving local coordinate charts, and by definition of integration of volume forms, the measure $d(kT)$ is absolutely continuous with respect to the Lebesgue measure on $\mathbb{R}^{\mathrm{dim}(K/T)}$ when restricted to one of these charts $U$, see Chapter 6 in \cite{rudin} for the definition of absolute continuity. By the property of the Lebesgue measure, see Theorem 2.20 in \cite{rudin}, together with properties of a submanifold, see Section 1.1 in \cite{var}, we see that $\int_{E_j\cap U} h(kT)\, d(kT)=0$ for any of the submanifolds $E_j$ in $E$, and for any smooth function $h$ on $K/T$. Hence by collecting pieces, we get $\int_E h(kT)\, d(kT)=0$ for any smooth function $h$ on $K/T$, and thus also for any continuous function, again by denseness, and due to Lebesgue's dominated convergence theorem.    
\qed

Let's extract precisely what we need from this lemma, namely a preliminary integration formula which we will apply to compact connected simple Lie groups.

\begin{cor}
\label{corint}
Let $K$ be a compact Lie group with a closed subgroup $T$. Let $E\subset K/T$ be as in the preceding
Lemma. Assume that the complement $E^c$ in $K/T$ of $E$ is a submanifold of $K/T$, and let
$\iota\colon E^c\to K/T$ denote the inclusion map. Let $\omega$ be a normalized $K$-invariant analytic volume
form on $K/T$, and let $d(kT)$ be the Borel measure defined by $\omega$. Then the normalized Haar integral of
$f\in\mathbb{C}[K]$ is given by the formula  
$$
\int_K f(k)\, dk =\int_{E^c}\Phi^{-1}P(f)(kT)\, d(kT)=\int_{E^c}\iota^* (\Phi^{-1}P(f)\omega  )
=\int_{E^c}((\Phi^{-1}P(f)\iota )\, \iota^*(\omega ).
$$
\end{cor}


\section{Some Poisson geometry of Lie groups}\label{s-Poisson}
Let $K$ be any compact connected simple Lie group, and let $G$ be a complexification as provided by Theorem
\ref{th-complexif}. In this section we follow conventions and partly notation from \cite{lu} using results
from that article freely. 

Let $\mathfrak{g}$ be the Lie algebra of $G$. Choose a Cartan subalgebra $\mathfrak{h}$ of $\mathfrak{g}$ and
let $H$ be the corresponding connected Lie subgroup of $G$. Let $\Delta$ be the set of roots of $\mathfrak{g}$
with respect to $\mathfrak{h}$, and let $\Delta_+$ be the positive roots, writing $\beta>0$ for
$\beta\in\Delta_+$. Let $\mathfrak{b}$ be the Borel subalgebra of $\mathfrak{g}$ spanned by $\mathfrak{h}$ and
the positive root vectors, and let $B$ be the corresponding connected Lie subgroup of $G$.  

Let $(\cdot ,\cdot)$ be the Killing form on $\mathfrak{g}$. For $\alpha\in\Delta_+$, let $H_\alpha$ be the
image of $\alpha$ under the isomorphism $\mathfrak{h}^*\to\mathfrak{h}$ defined by the Killing form, so  
$(H_\alpha ,H)=\alpha (H)$ for all $H\in\mathfrak{h}$. Choose root vectors $E_{\pm\alpha}$ for $\pm\alpha$
such that $(E_\alpha ,E_{-\alpha})=1$, so $[E_\alpha ,E_{-\alpha}]=H_\alpha$. Then the real subspace 
$\mathfrak{k}$ of $\mathfrak{g}$ spanned over $\mathbb{R}$ by the vectors $iH_\alpha , E_\alpha-E_{-\alpha}$
and $i(E_\alpha +E_{-\alpha})$ for all $\alpha >0$ is a compact real form of $\mathfrak{g}$ with a compact
connected subgroup $K$ of $G$, and $T=K\cap B$ is a maximal compact torus of $K$ with Lie algebra
$\mathfrak{t}$ spanned over $\mathbb{R}$ by the vectors $iH_\alpha$ for all $\alpha >0$.  

Let $\mathfrak{a}=i\mathfrak{t}$ and let $\mathfrak{n}$ be the subalgebra of $\mathfrak{g}$ spanned over $\7R$
by the positive root vectors. Denote the corresponding connected Lie subgroups of $G$ by $A$ and $N$. Then as
real manifolds we have the Iwasawa decomposition $G=KAN$ saying that the product map 
$$
K\times A\times N\to G
$$
is a diffeomorphism. Note that $B=TAN=ANT$. The quotient map $G\to G/B$ restricted to $K\subset G$ induces an 
isomorphism $K/T\to G/B$ of the flag manifolds. In the sequel we may suppress this isomorphism.

Let $\alpha_i$ be the simple roots of $\mathfrak{g}$ with corresponding fundamental reflections 
$$
s_i\colon\mathfrak{h}^*\to\mathfrak{h}^*,\ \  x\mapsto x-2\frac{(x,\alpha_i )}{(\alpha_i ,\alpha_i )}\alpha_i ,
$$ 
where the Killing form has been dualized using $\mathfrak{h}^*\cong\mathfrak{h}$. The elements $s_i$ generate
the Weyl group $W$ of $\mathfrak{g}$. Its elements preserve the roots of $\mathfrak{g}$ and the dimensions of
the root spaces. The length $l(w)$ of a Weyl element $w\in W$ is the number of fundamental reflections
appearing in a reduced expression of $w$, that is, one with a minimal number of fundamental reflections. The
Weyl group is a finite Coxeter group. It has a longest element $w_0$ with length $l(w_0)$ equal to the number
of positive roots of $\mathfrak{g}$. The Weyl group can also be realized as the groups $N(H)/H$ and $N(T)/T$,
where $N(H)$ and $N(T)$ are the normalizers of $H$ and $T$ in $G$ and $K$, respectively. We write $\dot{w}$
for a representative of $w$ in any of the normalizers.

Let $w\in W$ and let $\Delta_+^w$ be the positive roots $\alpha$ of $\mathfrak{g}$ such that $w^{-1}(\alpha)$
are not positive. If $w=s_{i_1}\cdots s_{i_{l(w)}}$ is a reduced expression for $w$, then each element
$\beta_j =s_{i_1}\cdots s_{i_{j-1}}(\alpha_{i_j})$ for $j=1,\dots, l(w)$ occur exactly once in $\Delta_+^w$,
and they exhaust these positive roots. 
The complex span $\mathfrak{n}_w$ of the corresponding vectors $E_\alpha$ is a subalgebra of $\mathfrak{g}$
with connected Lie subgroup $N_w$ of $G$. We have $N_w =N\cap wN_- w^{-1}$, where $N_-$ is the `opposite' of
$N$, so it is the unipotent subgroup of $G$ which is dual to $N$ with respect to the complex torus $H=TA$. In
particular, we see that $N_w\subset N$. Also note that $\Delta_+^{w_0}=\Delta_+$ and that $N=N_{w_0}$. 

The group $G$ can be decomposed into double cosets indexed by the Weyl group, a partition 
$$
G=\bigcup_{w\in W} B\dot{w}B
$$ 
known as the Bruhat decomposition. This in turn gives rise to a decomposition
of the flag manifold $G/B$ into a finite disjoint union $\bigcup_{w\in W} B\dot{w}B /B$ of submanifolds, or
cells 
$$
\Sigma_w =B\dot{w}B/B
$$ 
known as Schubert cells. These $2l(w)$-dimensional cells form a CW-complex for
$G/B$ having a unique cell $\Sigma_{w_0}$ of greatest (real) dimension $2l(w_0)$.  
Note that we may write the Bruhat decomposition of $G$ as a disjoint union
$$
G=\bigcup_{w\in W,\, t\in T}AN\dot{w}ANt
$$ 
since $\dot{w}\in N(T)$.     

We get a left action 
$$
G\times K\to K\colon (g,k)\mapsto g\circ k
$$ 
of $G$ on $K$, by projecting onto the
compact part of $gk\in KAN$ in the Iwasawa decomposition of $G$. Then the $N_w$-orbit  
$C_{\dot{w}} =N_w\circ\dot{w}$ through $\dot{w}\in K$ is diffeomorphic to $\Sigma_w$ via the quotient map
$\pi\colon K\to K/T$. In fact, it is even a symplectic diffeomorphism between symplectic leaves. The
corresponding Poisson manifold for the symplectic leaves $C_{\dot{w}}$ is one for which it is a Poisson Lie
group, i.e. one that makes the product map a Poisson map, and it is (among infinitely many such Poisson
structures on $K$) known as the standard one. It has Poisson bracket 
$$
\{f,g\}=w_K (df,dg)
$$ 
for $f,g \in C^\infty (K)$ with bivector $(w_K )_k =r_r-r_l$, where 
$$
r=\frac{1}{2}\sum_{\alpha >0}
X_\alpha\wedge Y_\alpha
$$ 
is the standard classical $r$-matrix for $K$, while $r_r$ and $r_l$ are respectively
the right- and left invariant bivector fields on $K$ with values $r$ at the unit element of $K$. Since 
the Poisson bivector $w_K$ vanishes on $T$, the homogeneous space $K/T$ has a unique Poisson structure such
that the quotient map $\pi$ and the action $K\times K/T\to K/T$ are Poisson maps. Any Poisson manifold is a
disjoint union of (immersed) symplectic submanifolds, each being the maximal one for which any two points can
be joined by a piecewise smooth curve, where each smooth segment is part of an integral curve of a Hamiltonian
vector field. Now $\Sigma_w$ is one such leaf in the Poisson flag manifold $K/T$.

While for $K/T$ there is only one leaf for each $w\in W$, this is not so for $K$. We have exactly one leaf
$C_{\dot{w}}$ that goes through every representative $\dot{w}\in N(T)\subset K$ of $w\in W=N(T)/T$. Let
$t\in T$. Then the leaf in $K$ that goes through the representative $\dot{w}t$ of $w$ is actually
$C_{\dot{w}}t$, and obviously $\pi (C_{\dot{w}}t)=\pi (C_{\dot{w}})=\Sigma_w$. The Bruhat piece
$K_w=K\cap (B\dot{w}B)$ is a disjoint union 
$$
\bigcup_{t\in T} C_{\dot{w}}t
$$ 
of leaves, and $K$ is a disjoint
union $\bigcup_{w\in W}K_w$ of such pieces, showing that the leaves of $K$ are parametrized by the set
$W\times T$.

In fact, the leaves of $K$ are exactly the $AN$-orbits under the previous action $\circ$ of $G$ on $K$. When
restricted to $AN$ this action is known as the dressing action on $K$ of the Poisson dual $K^*$ of $K$ since
$K^*$ and $AN$ are isomorphic as Poisson Lie groups. Let $w\in W=N(T)/T$. It happens to be the case that the
$AN$-orbit of $\dot{w}\in K$ is $N_w\circ\dot{w}=C_{\dot{w}}$. As an $AN$-homogeneous space $N_w$ is
diffeomorphic to $AN/AN^+_w$ thanks to the decomposition $AN=N_w AN^+_w$, where $N_w^+ =N\cap wNw^{-1}$. From
this and the Bruhat decomposition above, we see that $C_{\dot{w}}\subset AN\dot{w}AN$ and that every
$AN$-orbit is of the form 
$$
AN\dot{w}ant
$$ 
for $a\in A$, $n\in N$ and $t\in T$. So the leaf in $K$ that
goes through $\dot{w}t$ is indeed $C_{\dot{w}}t$.   

Consider the diffeomorphism $J_w\colon N_w\to C_{\dot{w}}$ given by 
$$
J_w (n)=n\circ\dot{w},
$$ 
and the
(holomorphic) diffeomorphism $j_w\colon N_w\to \Sigma_w$ given by 
$$
\iota_w j_w =\pi J_w ,
$$ 
where $\iota_w\colon\Sigma_w\to K/T$ is the inclusion map. Let $n\in N_w$. Using the Iwasawa decomposition, we may write $n\dot{w}=kb$ for $k\in K$ and $b\in AN$, and 
then $n\circ\dot{w}=k$. Hence $\pi J_w (n)=\pi (n\circ\dot{w})=\pi (k)=kB=kbB=n\dot{w}B$. 

Let $w\in W$ with a representative $\dot{w}\in N(T)\subset K$. For $n\in N_w$ let $a_w (n)$ be the $A$-component in the Iwasawa decomposition of $\dot{w}^{-1}n\dot{w}\in G=KAN$. This gives a differential map 
$$
a_w\colon N_w\to A ,
$$  
which is independent of the representative of $w$. It can also be pushed down to a well-defined differential map $C_{\dot{w}}\to A$ by $a_w (n\circ\dot{w})=a_w (n)$ for $n\in N_w$.        

Let $\lambda\in\mathfrak{h}^*$ be a weight of a finite dimensional representation of $\mathfrak{g}$.
Whenever $x\in A$, we can define a complex number $x^\lambda$ as follows:
Note that the exponential map $\mathrm{exp}\colon\mathfrak{a}\to A$ is surjective since $A$ is connected and Abelian, see e.g.\
\cite[Remark on p. 88]{var}. So we may write $x=\mathrm{exp}(X)$ for some $X\in\mathfrak{a}\subset\mathfrak{h}$. Then set $x^\lambda =e^{\lambda (X)}$. This is consistent with the formal manipulation $x^\lambda =e^{\lambda\, \mathrm{ln}x}=e^{\lambda\, \mathrm{ln}\, \mathrm{exp}(X)}=e^{\lambda X}$. The number $e^{\lambda (X)}$ is independent of what $X\in\mathfrak{a}$ we picked thanks to Theorem 5.107, Proposition 5.4 and Proposition 4.58 in \cite{knapp}.  

Let $dn$ denote the Haar measure on $N_w$ which we normalize such that 
$\int_{N_w}a_w (n)^{-4\rho}\, dn =1$, where $\rho$ is half times the sum of the positive roots of $\mathfrak{g}$, and $\rho$ is indeed a weight, even a dominant weight, see \cite[Proposition 21.16]{bump}. 

For a simple root $\alpha_i$ of $\mathfrak{g}$, let $c_{\alpha_i} =2/(\alpha_i ,\alpha_i )$, and note that 
\[ \begin{pmatrix} 1 & 0 \\ 0 & -1 \end{pmatrix} \mapsto c_{\alpha_i}\, H_{\alpha_i},\quad
\begin{pmatrix} 0 & 1 \\ 0 & 0 \end{pmatrix} \mapsto \sqrt{c_{\alpha_i}}\,  E_{\alpha_i},\quad
\begin{pmatrix} 0 & 0 \\ 1 & 0 \end{pmatrix} \mapsto \sqrt{c_{\alpha_i}}\,  E_{-\alpha_i}
\]
extends by $\mathbb{C}$-linearity to a Lie algebra homomorphism $\mathfrak{sl}(2,\mathbb{C})\to\mathfrak{g}$,
which induces a Lie group homomorphism 
\be \varphi_{\alpha_i}\colon SL(2,\mathbb{C})\to G \label{eq-varphi}\ee
such that $\varphi_{\alpha_i} (SU(2))\subset K$.

We discuss briefly the image of this homomorphism. Consider any parabolic subgroup $P$ of $G$, that is, a
proper subgroup of $G$ containing $B$. Then $K\cap P$ is a Poisson-Lie subgroup of $K$. Indeed, it suffices to
check that it is invariant under the dressing action, and to this end, let $k\in K\cap P$ and $b\in AN$.  
By using the Iwasawa decomposition, we may write $bk=k_1b_1$ for $k_1\in K$ and $b_1\in AN$. Then
$b\circ k=k_1 =bkb_1^{-1}\in K\cap P$ as $k_1\in K$ and $b,b_1\in AN\subset B\subset P$.   

There is a Poisson structure on $K/(K\cap P)\cong G/P$ rendering the quotient map $K\to K/(K\cap P)$
and the action $K\times K/(K\cap P)\to K/(K\cap P)$ Poisson maps. Moreover, 
the decomposition into symplectic leaves coincides with the Schubert cell decomposition  
$$
G/P =\bigcup_{w\in W/W_P}AN\dot{w} P/P 
$$ 
of $G/P$, where $W_P$ is the Weyl group of the Levi factor of $P$. This result generalizes the previous case where $P=B$. In another direction, let $\mathfrak{g}_{-\alpha_i}$ be the root space of $-\alpha_i$. Let $P_i$ be the connected subgroup of $G$ with Lie algebra $\mathfrak{b}\oplus \mathfrak{g}_{-\alpha_i}$. Clearly $P_i$ is a parabolic subgroup of $G$, at least when $G$ is bigger than $SL(2,\mathbb{C})$, and $K\cap P_i =\varphi_{\alpha_i} (SU(2))$. Now $SU(2)$ is a Poisson-Lie group with bivector $w_{SU(2)}$ defined from the $r$-matrix 
$$
\frac{1}{4}(\alpha_i ,\alpha_i )\begin{pmatrix}0&1\\-1&0\end{pmatrix}\wedge
\begin{pmatrix}0&i\\i&0\end{pmatrix}\in\mathfrak{su}(2)\wedge\mathfrak{su}(2). 
$$    
Identifying $\mathfrak{a}+\mathfrak{n}$ with $\mathfrak{k}^*$ by using the imaginary part of the Killing form, one sees that the subspace of $\mathfrak{a}+\mathfrak{n}$ that annihilates $d\varphi_{\alpha_i} (\mathfrak{su}(2))$ is an ideal of $\mathfrak{a}+\mathfrak{n}$. Hence $\varphi_{\alpha_i}$ is an isomorphism between the Poisson-Lie groups $SU(2)$ and $K\cap P_i$.  

\bigskip

Finding suitable coordinates for $C_{\dot{w}}$ and $\Sigma_w$ requires a closer look at $N_w$. Consider a reduced expression $w=s_{i_1}\cdots s_{i_{l(w)}}$ for $w\in W$, where $s_{i_j}$ is the fundamental reflection associated to a simple root $\alpha_{i_j}$ of $\mathfrak{g}$. Now 
$\begin{pmatrix}
0&i\\
i&0
\end{pmatrix}$ is a representative of the non-trivial element in the Weyl group of $SL (2,\mathbb{C})$. Therefore
$$\dot{s}_{i_j} =\varphi_{\alpha_{i_j}} \left(\begin{pmatrix}
0&i\\
i&0
\end{pmatrix} \right)\in K 
$$
is a representative of $s_{i_j}$, and hence 
$\dot{w}=\dot{s}_{i_1}\dot{s}_{i_2}\cdots \dot{s}_{i_{l(w)}}\in K$ is a representative of $w$. We have the following beautiful decomposition result, see \cite[Theorem 1]{lu}:

\begin{prop}
There exists a diffeomorphism 
\[ F_w\colon N_{s_{i_1}}\times N_{s_{i_2}}\times\cdots\times N_{s_{i_{l(w)}}}\to N_w \]
between these real manifolds. This map is characterized by the identity
$$
F_w (n_1,n_2,\dots,n_{l(w)})\circ w =(n_1\circ\dot{s}_{i_1})(n_2\circ\dot{s}_{i_2})\cdots (n_{l(w)}\circ\dot{s}_{i_{l(w)}})\in K,
$$
which says that $C_{\dot{w}}=N_w\circ\dot{w}$ is diffeomorphic to
$C_{\dot{s}_{i_1}}C_{\dot{s}_{i_2}}\cdots C_{\dot{s}_{i_{l(w)}}}$, where the product between the symplectic
leaves is the usual one in $K$. In fact, this latter diffeomorphism is even a symplectic diffeomorphism, but
it is in general not holomorphic. 
\end{prop}

We are now in a position to produce coordinates:

\begin{cor}
Keeping the notation, the map 
$$
\psi_w\colon\mathbb{C}^{l(w)}\to N_w,\ \ (z_1, z_2,\dots ,z_{l(w)})\mapsto F_w (n_{z_1}, n_{z_2},\dots, n_{z_{l(w)}})
$$  
gives coordinates $\{z_1 ,\bar{z}_1, z_2 , \bar{z}_2 , \dots , z_{l(w)}, \bar{z}_{l(w)}\}$ of $N_w$ as a real manifold, where 
$$
n_{z_j} =\mathrm{exp}\, (z_j \sqrt{c_{\alpha_{i_j}}}\,  E_{\alpha_{i_j}})
=\varphi_{\alpha_{i_j}}\left (\begin{pmatrix}
1&z_j\\
0&1
\end{pmatrix} \right)\in N_{s_{i_j}}
$$
for $j=1,\dots ,l(w)$.
\end{cor}

Note that these coordinates do depend on the reduced expression for $w$ since $F_w$ does. To see what these coordinates look like on the level of the leaves in $K$, write $z_j =z$ and $x=(1+|z|^2 )^{-1/2}$ for simplicity. Then the Iwasawa decomposition
$$
\begin{pmatrix}
1&z\\
0&1
\end{pmatrix}
\begin{pmatrix}
0&i\\
i&0
\end{pmatrix}
=\begin{pmatrix}
izx& ix\\
ix&-i\bar{z}x
\end{pmatrix}
\begin{pmatrix}
1/x &\bar{z}x\\
0&x
\end{pmatrix}
$$ 
in $SL(2,\mathbb{C})$ gives the following Iwasawa decomposition
$$
n_z\dot{s}_{i_j} =\varphi_{\alpha_{i_j}}\left( \begin{pmatrix}
izx& ix\\
ix&-i\bar{z}x
\end{pmatrix} \right)\varphi_{\alpha_{i_j}} \left( \begin{pmatrix}
1/x&\bar{z}x\\
0&x
\end{pmatrix} \right)
$$ 
in $G$. Thus the map
$$
\mathbb{C}\to C_{\dot{s}_{i_j}}, \ \ z\mapsto n_z\circ\dot{s}_{i_j}
=\varphi_{\alpha_{i_j}} \left( \begin{pmatrix}    izx& ix\\  ix&-i\bar{z}x  \end{pmatrix} \right)
=\varphi_{\alpha_{i_j}} \left( \frac{i}{(1+|z|^2))^{1/2}}\begin{pmatrix} z & 1 \\ 1 & -\bar{z}\end{pmatrix}\right)
$$
is a parametrization of the factor $C_{\dot{s}_{i_j}}$ in $C_{\dot{w}}$ in terms of the coordinates
$\{z,\bar{z}\}=\{z_j,\bar{z}_j\}$. 
 
Recall that $\beta_j =s_{i_1}\cdots s_{i_{j-1}}(\alpha_{i_j})$ for $j=1,\dots, l(w)$, when
$w=s_{i_1}\cdots s_{i_{l(w)}}$ is a reduced expression for $w\in W$. Then each 
$\gamma_j =-w^{-1}\beta_j = s_{i_{l(w)}}s_{i_{l(w)-1}}\cdots s_{i_{j+1}}(\alpha_{i_j})$ for $j=1,\dots, l(w)$,
occur exactly once in $\Delta_+^{w^{-1}}$, and they exhaust these positive roots.   

Denote by $\Omega_w$ the symplectic 2-form on $\Sigma_w$, and use the symbol $\mu_w$ for the Liouville volume
form $\frac{1}{l(w)!}\Omega_w^{l(w)}$. In terms of the coordinates
$\{z_1 ,\bar{z}_1, z_2 , \bar{z}_2 , \dots , z_{l(w)}, \bar{z}_{l(w)}\}$ we have the following expressions,
see Theorems 2, 3, 9 and Corollary 11 in \cite{lu}: 

\begin{thm}
 Keeping notation as above, we have
 $$
 \psi_w^* j_w^* (\Omega_w ) (z_1,\dots ,z_{l(w)})=\sum_{j=1}^{l(w)}\frac{i}{(\beta_j,\beta_j )}
 \frac{1}{1+|z_j|^2}\, dz_j\wedge d\bar{z}_j
 $$
 $$
 \psi_w^* j_w^* (\mu_w )(z_1,\dots ,z_{l(w)}) =\prod_{j=1}^{l(w)}\frac{i}{(\beta_j,\beta_j )}
 \frac{1}{1+|z_j|^2}\, dz_j\wedge d\bar{z}_j
$$
$$
\psi_w^* (dn)(z_1,\dots ,z_{l(w)}) =\lambda_w\prod_{j=1}^{l(w)}(1+|z_j|^2)
^{\frac{2(\rho ,\gamma_j )}{(\gamma_j,\gamma_j )}-1}
 \, dz_j\wedge d\bar{z}_j ,
$$ 
where
$$
\lambda_w =(i/\pi)^{l(w)}\prod_{j=1}^{l(w)}\frac{(\rho ,\gamma_j )}{(\gamma_j,\gamma_j )}
$$
and we think of the Haar measure $dn$ on $N_w$ as a volume form. It is normalized such that 
$$
\int_{N_w}a_w (n)^{-4\rho}\, dn =1.
$$
Furthermore, we have
$$
\psi_w^* (a_w )(z_1,\dots ,z_{l(w)}) =\prod_{j=1}^{l(w)}\exp\Big( \frac{1}{2}\mathrm{ln}(1+|z_j|^2 )c_j H_{\gamma_j} \Big),
$$
where $c_j=c_{\alpha_{i_j}}=2/(\alpha_{i_j},\alpha_{i_j})$, and the intrinsic relation
$$
j_w^* (\mu_w )=\Big(\prod_{j=1}^{l(w)}\frac{\pi}{(\rho ,\gamma_j )} \Big)\, a_w (n)^{-2\rho}\, dn .
$$
\end{thm} 

The Kostant harmonic differential forms $s^w$ are certain closed $K$-invariant analytic (as analytic translates of elements at $eT$) forms on the flag manifold $K/T$ of degree $2l(w)$. Their cohomology classes as $w$ ranges over the Weyl group $W$, furnish a basis for the de Rham cohomology $H^* (K/T,\mathbb{C})$ of $K/T$, which up to scalars is dual to the basis for the homology $H_* (K/T,\mathbb{C})$ with homology classes consisting of the closures of the Schubert cells $\Sigma_w\subset K/T$. The scalars are given as follows
$$
\int_{\Sigma_w} \iota_w^* (s^w )=\prod_{j=1}^{l(w)}\frac{\pi}{(\rho ,\beta_j )}.
$$   
We introduce a normalization constant $c_w$ which is set to be the inverse of this scalar.
We won't explain how the Kostant harmonic forms are constructed, nor discuss their nice properties, but we
will describe them in terms of our coordinates, see \cite[Theorem 14]{lu}, since this will be useful. 
  
\begin{thm}
 With notation as above, we have
 $$
\psi_w^*j_w^* \iota_w^* (s^w)(z_1,\dots ,z_{l(w)}) =\prod_{j=1}^{l(w)}\frac{i}{(\beta_j,\beta_j )}(1+|z_j|^2)
^{-\frac{2(\rho ,\beta_j )}{(\beta_j,\beta_j )}-1}
 \, dz_j\wedge d\bar{z}_j .
 $$
 Moreover, we have the intrinsic relation
 $$
j_w^*\iota_w^* (s^w ) =(a_w)^{2w^{-1}\rho} j_w^* (\mu_w ).
 $$
\end{thm}

In particular, since for the longest element $w_0$ of the Weyl group, we have $w_0^{-1}\rho =w_0\rho =-\rho$
and $\Delta_+^{w_0} =\Delta_+$, and $2l(w_0 )=\mathrm{dim}(K/T)$, we get the following result: 
 
 \begin{cor}
 \label{corkos}
With notation as above, the Kostant harmonic form $s^{w_0}$ is a $K$-invariant analytic volume form on $K/T$, and we have
$$
\psi_{w_0}^*j_{w_0}^* \iota_{w_0}^* (s^{w_0})(z_1,\dots ,z_{l(w_0)}) =\prod_{j=1}^{l(w_0 )}\frac{i}{(\beta_j,\beta_j )}(1+|z_j|^2)
^{-\frac{2(\rho ,\beta_j )}{(\beta_j,\beta_j )}-1}
 \, dz_j\wedge d\bar{z}_j 
 $$  
 with intrinsic relations
 $$
 j_{w_0}^*\iota_{w_0}^* (s^{w_0} )=(a_{w_0})^{-2\rho}\, j_{w_0}^* (\mu_{w_0})=
 \Big(\prod_{\alpha >0}\frac{\pi}{(\rho, \alpha )}\Big)\, a_{w_0}(n)^{-4\rho}\, dn .
 $$
\end{cor}

The corollary is consistent with the formula (5.8) for the Haar measure on $K$ given in \cite{resh-yak}, except that the normalization constant there is  erroneously inverted.


\section{Integration formulae for compact connected Lie groups}\label{sec-int}

We now combine the results obtained so far to obtain an explicit formula for the Haar integral on a compact
connected Lie group $K$. We proceed in a few steps. 

Consider first a compact connected simple Lie group $K$ with a maximal  torus $T$. Invoke the $K$-invariant
analytic volume form $\omega =c_{w_0}s^{w_0}$ on $K/T$, and consider the union $E=\bigcup_{w\neq w_0}\Sigma_w$
of submanifolds of $K/T$ all of dimension less than $\mathrm{dim}(K/T)$, and note that $E^c
=\Sigma_{w_0}$. Plugging this into the formula in Corollary \ref{corint}, and invoking Corollary \ref{corkos},
we get for $f\in\mathbb{C}[K]$ that 
$$
c_{w_0}^{-1}\int_K f(k)\, dk =\int_{\Sigma_{w_0}}((\Phi^{-1}P(f)\iota_{w_0})\, \iota_{w_0}^*(s^{w_0})
=\int_{\mathbb{C}^{l(w_0 )}}\psi_{w_0}^*j_{w_0}^*((\Phi^{-1}P(f)\iota_{w_0})\, \psi_{w_0}^*j_{w_0}^* \iota_{w_0}^* (s^{w_0}).
$$
Now $\iota_{w_0} j_{w_0}=\pi J_{w_0}$, so $j_{w_0}^*\iota_{w_0}^* =J_{w_0}^*\pi^*$ and $\Phi^{-1}(g)=(\pi^* )^{-1} (g)$ for $g$ analytic in $C^T (K)$. Hence
\begin{eqnarray}
  \lefteqn{ c_{w_0}^{-1}\int_K f(k)\, dk = \int_{\mathbb{C}^{l(w_0 )}}\psi_{w_0}^*J_{w_0}^*(P(f))\, \psi_{w_0}^*j_{w_0}^* \iota_{w_0}^* (s^{w_0}) }
   \nn\\
  &=& \int_{\mathbb{C}^{l(w_0 )}}P(f) ((n_{z_1}\circ \dot{s}_{i_1})(n_{z_2}\circ \dot{s}_{i_2})\cdots
(n_{z_{l(w_0 )}}\circ \dot{s}_{i_{l(w_0)}}))\, \psi_{w_0}^*j_{w_0}^* \iota_{w_0}^* (s^{w_0})
(z_1 ,z_2 ,\dots ,z_{l(w_0)}) \nn\\
  &=& \int_{\mathbb{C}^{l(w_0 )}}P(f) ((n_{z_1}\circ \dot{s}_{i_1})(n_{z_2}\circ \dot{s}_{i_2})\cdots
(n_{z_{l(w_0 )}}\circ \dot{s}_{i_{l(w_0)}}))\, \prod_{j=1}^{l(w_0 )}\frac{i}{(\beta_j,\beta_j )}(1+|z_j|^2)
^{-\frac{2(\rho ,\beta_j )}{(\beta_j,\beta_j )}-1}
 \, dz_j\wedge d\bar{z}_j \nn\\
  &=& \int_{\mathbb{C}^{l(w_0 )}}\int_T f((n_{z_1}\circ \dot{s}_{i_1})(n_{z_2}\circ \dot{s}_{i_2})\cdots
(n_{z_{l(w_0 )}}\circ \dot{s}_{i_{l(w_0)}})t)\, dt\, \prod_{j=1}^{l(w_0 )}\frac{i}{(\beta_j,\beta_j )}(1+|z_j|^2)
^{-\frac{2(\rho ,\beta_j )}{(\beta_j,\beta_j )}-1}
      \, dz_j\wedge d\bar{z}_j  \nn\\
  &=& \int_{\mathbb{C}^{l(w_0 )}}\int_T f\Big(\Big[\prod_{j=1}^{l(w_0)}\varphi_{\alpha_{i_j}}\bigg( \frac{i}{(1+|z_j|^2)^{1/2}}
      \begin{pmatrix} z_j & 1\\ 1 & -\bar{z}_j\end{pmatrix}\bigg)\Big] t\Big) \, dt\, \prod_{j=1}^{l(w_0 )}\frac{i}{(\beta_j,\beta_j )}(1+|z_j|^2)
^{-\frac{2(\rho ,\beta_j )}{(\beta_j,\beta_j )}-1}
      \, dz_j\wedge d\bar{z}_j,  \nn \\ 
\label{eq-poisson}  
\end{eqnarray}
where the notation $\prod_{j=1}^J k_j$ stands for $k_1k_2\cdots k_J$.

\begin{rema} \label{rem1}
1. The final row (\ref{eq-poisson}) of the above computation is the conclusion of the Poisson geometric approach;
the rest of the paper is devoted to its fairly elementary  analysis. From now on we abbreviate, writing $L$
and $\varphi_j$ instead of $l(w_0)$ and $\varphi_{\alpha_{i_j}}$, respectively.

2. Eq.\ (\ref{eq-poisson}) involves an integration over the maximal torus $T$. In order to settle the
terminology, we briefly recall the salient facts. Given a compact torus $T$ with Lie algebra $\mathfrak{t}$,
the exponential map $\mathfrak{t}\to T$ is a surjective homomorphism from the additive group $\mathfrak{t}$ to
the multiplicative group $T$, and it has a discrete kernel $\Lambda$, so that with
$\mathbb{T}\equiv\{z\in\mathbb{C}\, |\, |z|=1\}$, and with $R=\dim(T)$ we get
\[ T\cong\mathfrak{t}/\Lambda\cong (\mathbb{R}/\mathbb{Z})^R\cong\mathbb{T}^R \]
as Lie groups, e.g.\ see Proposition 15.3 in \cite{bump}. In terms of the coordinates
$z=(z_1,\ldots,z_R)\in\mathbb{T}^R$, the normalized Haar measure is well known to be
$\frac{1}{(2\pi i)^R}\frac{dz_1}{z_1}\cdots \frac{dz_R}{z_R}$. The group $\widehat{T}$ of characters on $T$,
that is, the homomorphisms from $T$ to $\mathbb{T}$ (with group product given by pointwise multiplication),
are just the finite dimensional unitary irreducible representations of $T$. It can be identified with
$\mathbb{Z}^R$, where ${\bf n}=(n_1 ,\dots ,n_R)\in\mathbb{Z}^R$ corresponds to the character 
$z\mapsto z^{\bf{n}}\equiv\prod_{i=1}^R z_i^{n_i}$.   
Matrix coefficients of representations of $T$ will therefore, under these identifications, be Laurent
polynomials in the variables $z_i$. 
\end{rema}

Let $h\in C(SU(2))$ and consider
$$
I\equiv\frac{i}{(\beta ,\beta )}\int_{\mathbb{C}}
h \left( \frac{i}{(1+|z|^2)^{1/2}}\begin{pmatrix} z & 1 \\ 1 & -\bar{z} \end{pmatrix} \right) (1+|z|^2)
^{-\frac{2(\rho ,\beta )}{(\beta,\beta )}-1}
 \, dz\wedge d\bar{z}.
$$  
We change coordinates $(z,\bar{z})\to (x,\theta )$, where $x =(1+|z|^2 )^{-1/2}$ as before, and $e^{i\theta}=z/|z|$, where 
$x\in(0,1]$ and $\theta\in [0,2\pi)$. Then $x^2 (1+|z|^2 )=1$, so 
$z=e^{i\theta}|z|=e^{i\theta}x^{-1}(1-x^2 )^{1/2}$. Hence
$$
\frac{\partial z}{\partial x}=-e^{i\theta}x^{-2}(1-x^2 )^{-1/2},\ \  \  \  \  \
\frac{\partial z}{\partial \theta}=ie^{i\theta}x^{-1}(1-x^2 )^{1/2},
$$
so
$$
dz\wedge d\bar{z} = (\frac{\partial z}{\partial x} dx+\frac{\partial z}{\partial\theta}d\theta )\wedge
(\frac{\partial\bar{z}}{\partial x} dx+\frac{\partial\bar{z}}{\partial\theta}d\theta )
=(\frac{\partial\bar{z}}{\partial x}\frac{\partial z}{\partial\theta}-\frac{\partial z}{\partial x}\frac{\partial\bar{z}}{\partial\theta})\, d\theta\wedge dx 
=-2ix^{-3}\, d\theta\wedge dx .
$$
Plugging this in, we get
\bea\label{eqL}
I &=& \frac{2}{(\beta ,\beta )}\int_0^1\int_0^{2\pi}
h\left( \begin{pmatrix}  ie^{i\theta}\sqrt{1-x^2} & ix \\ ix & -ie^{-i\theta}\sqrt{1-x^2} \end{pmatrix} \right)
(x^2)^{\frac{2(\rho ,\beta )}{(\beta,\beta )}-1}  x\, d\theta dx \\
  &=& \frac{2}{(\beta ,\beta )} \int_0^1\int_{S^1}
h\left( \begin{pmatrix}  iw\sqrt{1-x^2} & ix \\ ix & -iw^{-1}\sqrt{1-x^2} \end{pmatrix} \right)
(x^2)^{\frac{2(\rho ,\beta )}{(\beta,\beta )}-1}x  \frac{dw}{iw}  dx
\eea

(The first formula is more suited for computations, the second for theoretical purposes.)

Applying this coordinate transformation to all $L$ variables $z_j$ in (\ref{eq-poisson}) we arrive at
\begin{eqnarray} \int_K f(k)\, dk &=& c_{w_0}\int_0^1\int_{S^1}\cdots\int_0^1\int_{S^1} \int_T 
   f\left(\bigg[\prod_{j=1}^L\varphi_j\left(
   \begin{pmatrix}iw_j\sqrt{1-x_j^2} & ix_j \\ ix_j & -iw_j^{-1}\sqrt{1-x_j^2} \end{pmatrix}\right) \bigg]t\right) dt  \nn\\
  && \cdot\, \prod_{j=1}^L(\delta_j (x_j)\, \frac{dw_j}{iw_j} dx_j ), \label{eq-poisson2}
\end{eqnarray}
where
\[ \delta_j (x_j )=\frac{2}{(\beta_j ,\beta_j )}(x_j^2) ^{\frac{2(\rho ,\beta_j )}{(\beta_j ,\beta_j )}-1} x_j.\]

Since every $f\in \7C[K]$ is integrable over $K$, the factors $\delta_j(x_j)$ must be integrable. This can be seen directly:

\begin{lem}\label{lem-exp}
For any positive root $\alpha$ of a simple complex Lie algebra $\mathfrak{g}$, we have
$$
\frac{2(\rho ,\alpha )}{(\alpha ,\alpha )}\in\mathbb{N}=\{1,2,\ldots\}, 
$$
where $\rho$ is half times the sum of all the positive roots, and $(\cdot ,\cdot )$ is the Killing form, or
any other invariant bilinear form. 
Thus each $\delta_j(x_j)$ in (\ref{eq-poisson2})  is an odd natural power of $x_j$.
\end{lem}

\bp
We know that $\rho$ is a dominant weight, and for weights the scalar in question is always an integer, see the
paragraph before Proposition 21.14 in \cite{bump}. Since $(\alpha ,\alpha )>0$,  and since $\alpha$ is a
finite sum $\sum n_i\alpha_i$ of simple roots $\alpha_i$, where $n_i$ are non-negative integers, it suffices
to show that each $(\rho ,\alpha_i )>0$. Now this is clear from the identities
$2(\rho ,\alpha_i )/(\alpha_i ,\alpha_i ) =1$, see the proof of Proposition 21.16 in \cite{bump}.  
\qed

\begin{rema}\label{rem2}
1. With $D=\dim(K)$ and the rank $R=\dim(T)$, in (\ref{eq-poisson2}) we have $(D-R)/2$ $[0,1]$-valued
variables $x_j$, equally many $S^1$-valued variables $w_j$, and $R$ $S^1$-valued variables $z_k$ parametrizing
$T$. (That $D-R$ is even already follows from the classical fact that $G/B\cong K/T$ admits a complex
structure.) Thus in (\ref{eq-poisson2}) we are integrating over $[0,1]^{(D-R)/2}\times \7T^{(D+R)/2}$.

2. It is instructive to consider the case $K=SU(2)$. The standard choice of the maximal torus is
$T=\{\mathrm{diag}(w,w^{-1})\ | \ w\in\7T\}$.
Note that we have only one simple root $\alpha$, so $2\rho =\alpha$ and $c_{w_0}=(\alpha ,\alpha )/2\pi$ and 
$\delta (x)=\frac{2}{(\alpha ,\alpha )}x$. Hence for $f\in\mathbb{C}[SU(2)]$, we get  
\[ \int_{SU(2)} f
  =\frac{1}{2\pi^2}\int_0^1\int_0^{2\pi}\int_0^{2\pi} f \left (\begin{pmatrix} ie^{i\theta}\sqrt{1-x^2} &ix\\ ix&   -ie^{-i\theta}\sqrt{1-x^2}\end{pmatrix}\begin{pmatrix} e^{i\eta} &0\\ 0& e^{-i\eta}\end{pmatrix} \right)
\,  d\eta d\theta\, xdx . \]
This formula coincides with the following classical one found, e.g., in  \cite[Ch.\ III, Sect.\ 6.1]{V}:
\be \int_{SU(2)} f =\frac{1}{16\pi^2}\int_0^{2\pi}\int_0^\pi\int_{-2\pi}^{2\pi}
  f\left( \begin{pmatrix} \cos\frac{\theta}{2}\,e^{\frac{i(\phi+\psi)}{2}} &
    i\sin\frac{\theta}{2}\,e^{\frac{i(\phi-\psi)}{2}} \\
    i\sin\frac{\theta}{2}\,e^{\frac{i(\psi-\phi)}{2}} & \cos\frac{\theta}{2}\,e^{\frac{-i(\phi+\psi)}{2}}
    \end{pmatrix}\right)\sin\theta \,d\psi\,d\theta\,d\phi. \label{eq-vil}\ee
To prove this equality, it suffices to consider the monomials
$\begin{pmatrix} a & c \\ b & d\end{pmatrix}\mapsto a^{n_1}b^{n_2}c^{n_3}d^{n_4}$, whose linear combinations
are uniformly dense in $C(SU(2))$. It is a simple matter to check that on these monomials both formulae evaluate to
$(-1)^{n_2}\delta_{n_1,n_4}\delta_{n_2,n_3}\int_0^1 x^{n_1}(1-x)^{n_2}dx$. For (\ref{eq-vil}) this computation
can be found in \cite{mueger-tus}, where the authors proved a reduction of the Mathieu conjecture for $SU(2)$
to a simpler conjecture free of Lie theory. The latter will be generalized to all compact connected Lie groups
in Section \ref{sec6}. 

3. Analogous integration formulae for other groups are hard to find in the mathematical literature. (Weyl's
celebrated integration formula goes in a different direction.) Traditionally, those interested in explicit
integration formulae have been theoretical physicists, primarily concerned with $SU(n)$, and proofs tend to
be sketchy. This lacuna was filled by the work of M.M.'s PhD student K.~Zwart \cite{zwart}, who proved
rigorous integration formulae for $G_2$ and for $SU(n), SO(n), Sp(n)$ and applied the results to the Mathieu
conjecture for these groups, finding sufficient conditions in Abelian terms. To produce parametrizations of
these groups, he used KAK decompositions. For the series $SU, SO, Sp$ (the spin groups unfortunately were not
discussed) this involves induction over $n$. While $G_2$ can be done `by hand', this is not feasible for the 
other exceptional groups due to their dimensions. Zwart's results involve the square roots that we eliminate,
see the comments following the proof of Lemma \ref{lem-1}. Thus his approach falls short of giving the
results of Section \ref{sec6}. 
\end{rema}

\begin{prop} Let $K$ be a compact connected simple Lie group and keep the other notations as before. Then the
formulae (\ref{eq-poisson}) and (\ref{eq-poisson2}) hold for all measurable $f:K\rarr\7C$. 
\end{prop}

\bp Formulae (\ref{eq-poisson}) and (\ref{eq-poisson2}) have been proven for all $f\in \7C[K]$. They extend to
$C(K)$ by uniform density of $\7C[K]\subset C(X)$ and to $L^1(K)$ by $L^1$-density of $C(X)\subset L^1(X)$.
\qed\\

Let now $K=K_1\times\cdots\times K_n\times T$, where the $K_i$ are compact connected simple Lie groups and $T$
is a compact torus. Since the normalized Haar measure on a direct product of compact groups is the product of
the normalized Haar measures, it is immediate how to integrate a function on $K$ by combining
(\ref{eq-poisson}) or (\ref{eq-poisson2}) with the considerations in Remark \ref{rem1}.2.

Finally, we can discuss integration of a function $f\in L^1(K)$ on an arbitrary compact connected Lie
group $K$. By Theorem \ref{thm-prod} there are compact connected simple Lie groups $K_1,\ldots, K_n$, a
compact torus $T$, a finite subgroup $D$ of the center of $K_1\times\cdots\times K_n\times T$ and an
isomorphism $K\rarr(K_1\times\cdots\times K_n\times T)/D$. Let
$p: K_1\times\cdots\times K_n\times T\rarr(K_1\times\cdots\times K_n\times T)/D$ be the quotient
map. Identifying $K$ and $(K_1\times\cdots\times K_n\times T)/D$, $f'=f\circ p$ is a measurable function on 
$K_1\times\cdots\times K_n\times T$ that we can integrate by the above methods. Finally we have 
\[ \int_K f=\int_{K_1\times\cdots\times K_n\times T}f' \]
by a special case of the results in Lemma \ref{lem-int}.


\section{Specialization to finite type functions}\label{sec5}
In this section we specialize to finite type functions.

\begin{lem} \label{lem-1}
Let $K$ be a compact connected simple Lie group with maximal torus $T$, and let $\psi: \7T^R\rarr K$ be a
homomorphism identifying $\7T^R$ with $T$. Let $w_0=s_{i_1}\cdots s_{i_L}$ be a reduced form of the maximal
length Weyl group element. Let $\varphi_j=\varphi_{\alpha_{i_j}}$ be the homomorphisms $SU(2)\to K$
constructed in Section \ref{s-Poisson} (around (\ref{eq-varphi})). Define
\bean  Q_K: &&  [0,1]^L\times \7T^L\times \7T^R\rarr K, \\
   &&  (x,w,z)\mapsto  \bigg[\prod_{j=1}^L\varphi_j\left(
    \begin{pmatrix}iw_j\sqrt{1-x_j^2} & ix_j \\ ix_j & -iw_j^{-1}\sqrt{1-x_j^2} \end{pmatrix}\right) \bigg]\psi(z)
\eean
Then for each $f\in \7C[K]$ the function $\widetilde{f}=f\circ Q_K: [0,1]^L\times \7T^L\times \7T^R\rarr\7C$
depends polynomially on the $x_j$ and $\sqrt{1-x^2_j}$ and is a Laurent polynomial in the $w_j$ and the  $z_k$, where
$j=1,\ldots,L,\, k=1,\ldots,R=\dim(T)$. The map $f\mapsto\widetilde{f}$ is a homomorphism of $\7C$-algebras. 
\end{lem}

\bp Since $f$ is finite type, there are a finite dimensional representation $(V,\Pi)$ of $K$ and $M\in\mathrm{End}\,V$ such that
$f(k)=\Tr_V(M\Pi(k))$ for all $k\in K$. Since $\Pi$ is a homomorphism, we have
\be \widetilde{f}(x,w,z)=  \Tr_V\left(M \bigg[\prod_{j=1}^L \big(\Pi\circ\varphi_j\big) \left(
   \begin{pmatrix}iw_j\sqrt{1-x_j^2} & ix_j \\ ix_j & -iw_j^{-1}\sqrt{1-x_j^2} \end{pmatrix}\right) \bigg](\Pi\circ\psi)(z)\right). \label{eq-ftilde}\ee
Now each $\Pi\circ\varphi_j$ is a representation of $SU(2)$ on the finite dimensional vector space $V$. It is
well known that each matrix element of $(\Pi\circ\varphi_j)(s)$ (with respect to any basis of $V$) is a
polynomial in the matrix elements of $s\in SU(2)$. (This is an immediate consequence of complete reducibility
and the fact that every irreducible representation of $SU(2)$ is contained in some tensor power of the
defining representation of $SU(2)$ on $\7C^2$.) Thus each factor $\Pi\circ\varphi_j(\cdots)$ is polynomial in
$x_j$ and $\sqrt{1-x^2_j}$ and Laurent w.r.t.\ $w_j$.
On the other hand, $\Pi\circ\psi$ is a representation of $\7T^R$ on $V$. Thus by the facts on the representations on tori mentioned in Remark \ref{rem1}.2, each matrix element of $(\Pi\circ\psi)(z)$, where $z\in\7T^R$, is a Laurent polynomial in $z_1,\ldots,z_R$.

Now the claim follows readily since the matrix elements of a product of matrices are polynomials in the entries of all the matrices.

The homomorphism property of $f\mapsto\widetilde{f}$ is evident.
\qed\\

Using the lemma, it is straightforward to reduce the Mathieu conjecture for compact connected simple Lie groups
and, in a second step, for all compact connected Lie groups to a Mathieu-style conjecture concerning functions that
are Laurent polynomials in a number of variables $z_i$ and polynomials in variables $x_j$ and
$\sqrt{1-x_j^2}$. Thanks to the specific form of (\ref{eq-poisson2}) one can, however, do better by removing 
the square roots. The following simple lemma is the key:

\begin{lem} \label{lem-2}
 Define
\bean   A: && SU(2)\rarr\7C^4,\ \  \begin{pmatrix} a & b \\ c & d \end{pmatrix}\mapsto(a,b,c,d) \\
   B: &&  [0,1]\times S^1\to\mathbb{C}^4, \ \ (x,w)\mapsto (iw(1-x^2 ), ix, ix, -iw^{-1}).
\eean
Then for any continuous function $\delta:(0,1]\to(0,\infty)$ and polynomial $p: \7C^4\rarr\7C$ we have
\[ \int_0^1\int_{S^1} pA\left(\begin{pmatrix}
iw\sqrt{1-x^2}& ix\\
ix&-iw^{-1}\sqrt{1-x^2}
\end{pmatrix}\right)\, \frac{dw}{w}\, \delta (x) dx
=\int_0^1\int_{S^1}pB (x,w)\, \frac{dw}{w}\, \delta (x) dx. \]
\end{lem}  

\bp
Since both sides are linear in $p$, it suffices to check the identity for 
$p=y_1^{n_1}\cdots y_4^{n_4}$, where the $n_i$ are non-negative integers. This amounts to proving
\[ \int_0^1\int_{S^1}w^{n_1-n_4}(1-x^2)^{(n_1+n_4)/2} x^{n_2+n_3}
\frac{dw}{w}\delta(x)\, dx
=\int_0^1\int_{S^1} w^{n_1-n_4}(1-x^2)^{n_1} x^{n_2+n_3}
\frac{dw}{w}\delta(x)\, dx, \]
which is obviously true since the integral over $S^1$ vanishes on both sides unless $n_1 =n_4$.
\qed


\begin{prop} \label{prop-tt} Let $K$ be a  compact connected simple Lie group, and keep the notations as above.
In particular $D$ and $R$ are the dimension and rank, respectively, of $K$ and $L=(D-R)/2$. Then there is a
linear map assigning to each $f\in \7C[K]$ a function $\stackrel{\approx}{f}$ that is polynomial in $L$
variables $x_j$ and Laurent polynomial in $L+R=(D+R)/2$ variables $w_j, z_k,$ and such that for all
$f,g\in \7C[K]$ and non-negative integers $a,b$, we have
\be \int_K f^ag^b = \frac{c_{w_0}}{(2\pi i)^R} \int_0^1\int_{S^1}\cdots\int_0^1\int_{S^1}\int_{\7T^R}
   \stackrel{\approx}{f}\!(x,w,z)^a\stackrel{\approx}{g}\!(x,w,z)^b
    \frac{dz_1}{z_1}\cdots\frac{dz_R}{z_R}\Big(\prod_{j=1}^L \frac{dw_j}{iw_j}\delta_j(x_j)dx_j \Big).
\label{eq-integ}\ee
\end{prop}

\bp In view of Lemma \ref{lem-1} and (\ref{eq-poisson2}) it is evident that (\ref{eq-integ}) is true for
$\widetilde{f}$ in place of $\stackrel{\approx}{f}$. 

It is clear from (\ref{eq-ftilde}) that $\widetilde{f}$ depends on $(x_j, w_j)$ only via the factor
$(\Pi\circ\varphi_j)(\cdots)$ and on $z$ only via $(\Pi\circ\psi)(z)$. Thus with (\ref{eq-poisson2}) we have
\be \int_K f= c_{w_0} \Tr_V\Big( M \Big[ \prod_{j=1}^L R_j \Big] P\Big), \label{eq-intf}\ee
where
\bean R_j &=& \int_0^1 \int_{S^1} (\Pi\circ\varphi_j)
    \left( \begin{pmatrix}iw_j\sqrt{1-x_j^2} & ix_j \\ ix_j & -iw_j^{-1}\sqrt{1-x_j^2} \end{pmatrix} \right) \, \frac{dw_j}{iw_j}\delta_j(x_j)dx_j, \\
  P &=& \frac{1}{(2\pi i)^R} \int_{\7T^R} (\Pi\circ\psi)(z) \frac{dz_1}{z_1}\cdots\frac{dz_R}{z_R}.
\eean
($P=\int_T \Pi(t)dt$ is an idempotent whose image is the subspace of $V$ of $\Pi(T)$-invariant vectors.)
Since $(\Pi\circ\varphi_j)\left(\begin{pmatrix}iw_j\sqrt{1-x_j^2} & ix_j \\ ix_j & -iw_j^{-1}\sqrt{1-x_j^2} \end{pmatrix} \right)$
is a matrix whose element $(k,l)$ is a polynomial $p_{j,k,l}A$ in the elements of the $2\times 2$ matrix, Lemma \ref{lem-1} gives
\[ (R_j)_{k,l}= \int_0^1\int_{S^1} p_{j,k,l}B(x_j,w_j) \frac{dw_j}{iw_j}\delta_j(x_j)dx_j. \]
Note that $p_{j,k,l}B(x_j,w_j)$ (and not only its integral) is independent of the choice of the polynomial
$p_{j,k,l}$: If $q$ is the difference of two such polynomials, $qA$ vanishes on $SU(2)$, thus also on the
complexification $SL(2,\7C)$, which contains $\begin{pmatrix} iw(1-x^2) & ix \\ ix &
-iw^{-1}\end{pmatrix}$ whose image under $A$ is $B(x,w)$.

Letting $C_j(x_j,w_j)$ be the matrix with elements $p_{j,k,l}B(x_j,w_j)$, we have
\be \stackrel{\approx}{f} (x,w,z)=\Tr_V\Big(M \Big[ \prod_{j=1}^L C_j(x_j,w_j) \Big] (\Pi\circ\psi)(z)\Big), \label{eq-ftt}\ee
which is polynomial in $x$ and Laurent in $w,z$ and independent of the chosen $p_{j,k,l}$.
By its construction and Lemma \ref{lem-2}, $\stackrel{\approx}{f}$ has the same integral as $\widetilde{f}$.

It is clear from the construction that $\stackrel{\approx}{cf}=c\stackrel{\approx}{f}$ for all $c\in\7C$. For
the proof of additivity of $f\mapsto\stackrel{\approx}{f}$, let
$f=\Tr_V(M\Pi(\cdot)),\,f'=\Tr_{V'}(M'\Pi'(\cdot))$. Replacing, if necessary, $(V,\Pi)$ and $(V',\Pi')$ by
$(V\oplus V',\Pi\oplus\Pi')$, we may assume $(V,\Pi)=(V',\Pi')$. Since $\stackrel{\approx}{f}$ is given by
(\ref{eq-ftt}) and $\stackrel{\approx}{f'}$ by the same formula with $M$ replaced by $M'$, it is clear that
$\stackrel{\approx}{(f+f')}=\stackrel{\approx}{f}+\stackrel{\approx}{f'}$, thus $f\mapsto\stackrel{\approx}{f}$
is linear.

Now let also $g\in \7C[K]$ and pick a representation $(V',\Pi')$ of $K$ and $M'\in\End V'$ such that
$g(k)=\Tr_{V'}(M'\Pi'(k))$. Applying the above considerations to
$fg=\Tr_{V\otimes V'}\big((M\otimes M')(\Pi\otimes\Pi')(\bullet)\big)$ gives
\bea \int_K fg &=& \int_K \Tr_{V\otimes V'}\big((M\otimes M')(\Pi\otimes\Pi')(k)\big)\,d\mu(k) \nn\\
  &=& c_{w_0} \Tr_{V\otimes V'}\Big((M\otimes M')\Big[\prod_{j=1}^L R''_j \Big] P''\Big),\label{eq-fg}
\eea
with
\bean
  P'' &=& \frac{1}{(2\pi i)^R} \int_{\7T^R} (\Pi\otimes\Pi')(\psi(z)) \frac{dz_1}{z_1}\cdots\frac{dz_R}{z_R}, \\
  R''_j &=& \int_0^1 \int_{S^1} (\Pi\otimes\Pi')\circ\varphi_j
    \left( \begin{pmatrix}iw_j\sqrt{1-x_j^2} & ix_j \\ ix_j & -iw_j^{-1}\sqrt{1-x_j^2} \end{pmatrix} \right) 
    \frac{dw_j}{iw_j}\delta_j(x_j)dx_j.
\eean
Since $\Pi\circ\varphi_j$ and $\Pi'\circ\varphi_j$ are finite dimensional
representations of $SU(2)$, there exist polynomials $p_{j,k,l}, p'_{j,k',l'}$ such that 
$(\Pi\circ\varphi_j)_{k,l}=p_{j,k,l}A$ and a similar primed identity hold. Applying Lemma
\ref{lem-2} to the matrix element $R''_{j,kk',ll'}$ gives
\[ R''_{j,kk',ll'}= \int_0^1 \int_{S^1} (p_{j,k,l}p'_{j,k',l'})(B(x_j,w_j)) \frac{dw_j}{iw_j}\delta_j(x_j)dx_j. \]
Due to the trivial fact $(p_{j,k,l}p'_{j,k',l'})(B(x_j,w_j))=p_{j,k,l}(B(x_j,w_j))p'_{j,k',l'}(B(x_j,w_j))$ we have
\[ R''_{j,kk',ll'}= \int_0^1 \int_{S^1} (C_j(x_j,w_j))_{k,l}(C'_j(x_j,w_j))_{k',l'} \frac{dw_j}{iw_j}\delta_j(x_j)dx_j, \]
or just
\[ R_j''= \int_0^1 \int_{S^1} C_j(x_j,w_j)\otimes C'_j(x_j,w_j) \frac{dw_j}{iw_j}\delta_j(x_j)dx_j. \]
Plugging this back into (\ref{eq-fg}) and comparing the result with (\ref{eq-ftt}) we obtain exactly the
r.h.s. of (\ref{eq-integ})  for $a=b=1$, completing the proof of this identity in this case. The general
case now follows by a simple induction argument.
\qed

\begin{rema}
1. We do not claim that $f\mapsto\stackrel{\approx}{f}$ is a homomorphism, but we expect that
$\stackrel{\approx}{(fg)}=\widetilde{f}\stackrel{\approx}{g}$ holds whenever $f,g\in \7C[K]$ and $\widetilde{f}$
is free of square roots.

2. One might wonder whether there is an upper bound to the powers of the variables in $\delta$ that is uniform
over the collection of compact connected simple Lie groups. But this is not the case since if we fix one
such group and consider its maximal positive root $\alpha$, we get the dual Coxeter number
$g=\frac{2(\rho,\alpha)}{(\alpha ,\alpha )} +1$, see the appendix in \cite{chari-press}, and according to
tables, see e.g.\ Wikipedia, this number equals $n+1$ for the series $A_n$, which includes $SU(n)$. So in this
case the corresponding variable in $\delta$ will occur with power $2n-1$, which has no upper bound as $n$
increases. 
\end{rema}

\vspace{.5cm}

We now generalize the above considerations to compact groups $K$ that are of the form
$K=K_1\times \cdots\times K_N\times T$, where the $K_n$ are compact connected simple Lie groups 
and $T$ is a torus. If $f_1,\ldots,f_N$ and $f_T$ are finite type functions on the groups $K_n$ and $T$, their 
product is a finite type function $f$ on $K$. Let $\psi_T:\7T^{\dim T}\rarr T$ be an isomorphism. Then it is
clear that 
\[ \stackrel{\approx}{f}(x,w,z,u)=\Big(\prod_{n=1}^N \stackrel{\approx}{f_n}(x_n,w_n,z_n)\Big) f_T(\psi_T(u)),
\]
where the $x_n,w_n,z_n,u$ clearly are multi-component variables, 
is polynomial in the $x$-variables and Laurent polynomial in the variables $w,z,u$.
By the preceding considerations, we have
\[ \int_K f=\big(\prod_{n=1}^n \int_{K_n}f_n\big)\int_T f_T=
   \big(\prod_{n=1}^n \int_{x_n,w_nz_n}\stackrel{\approx}{f_n}\big)\int_{\7T^{\dim T}} f_T\psi_T.\]
The notation is quite symbolic, suppressing the weight factors $\delta$ and the measures for the integrations
over $x_n,w_n,z_n,u$.) 

If $f'=f'_1\cdots f'_Nf'_T$ is another such product, we have
\[\int_K f^a{f'}^b=\big(\prod_{n=1}^N \int(\stackrel{\approx}{f_n})^a(\stackrel{\approx}{f'_n})^b\big)\int_T (f_T)^a(f'_T)^b.\]
Thus we have a map $f\mapsto\stackrel{\approx}{f}$ with the same properties as in Proposition \ref{prop-tt}, except
for the fact that so far it is defined only on the elements $f\in \7C[K]$ that are of the above product form. The
general element $f\in \7C[K]$ is a finite linear combinations of such products, and we define
$\stackrel{\approx}{f}$ by linearity. Now an easy argument involving the multinomial formula shows that
$\int_K f^a g^b=\int_{x,w,z,u} \stackrel{\approx}{f}^a\stackrel{\approx}{g}^b$ still holds.

We summarize the result, gathering variables of like type:

\begin{thm}\label{thm-main1}
Let $K$ be a direct product of finitely many compact connected simple Lie groups and a compact
torus. Then there are non-negative integers $N,M$, a monomial $\delta\in\7C[x_1,\ldots,x_N]$ of odd order in
each variable, a constant $C_K$ and a linear map
\[ \7C[K]\rarr \7C[x_1,\ldots,x_N,z_1^{\pm 1},\ldots,z_M^{\pm 1}],\ \ f\mapsto\stackrel{\approx}{f} \]
such that
\[ \int_K f^a g^b = C_K\int_{[0,1]^N\times\7T^M} \stackrel{\approx}{f}^a \stackrel{\approx}{g}^b\,
     \delta(x)\, dx_1\cdots dx_N \frac{dz_1}{z_1}\cdots\frac{dz_M}{z_M} \]
for all $f,g\in \7C[K]$ and non-negative integers $a,b$. In the sequel, we abbreviate the
r.h.s.\ to $C_K\int_{(N,M)} \stackrel{\approx}{f}^a \stackrel{\approx}{g}^b\! \delta$.
\end{thm}


\section{Reduction of the Mathieu conjecture}\label{sec6}
Theorem \ref{thm-main1} applies to all compact connected Lie groups that are direct products of simple groups
and circles. While in view of examples like $U(2)$ it is not the case that every compact connected Lie group is
of this form, Theorem \ref{thm-prod} shows that we are not far off. Indeed, for a discussion of the
Mathieu conjecture it suffices to consider groups of such product form:

\begin{prop}\label{prop-mathieu1}
If the Mathieu conjecture holds for some compact connected Lie group $K$, it holds for any quotient group
$K/N$, where $N$ is a closed normal subgroup of $K$. 
\end{prop}

\bp Let $p: K\rightarrow K/N$ be the quotient map. If $f:K/N\rightarrow\7C$ is a finite type function then
$\widehat{f}=f\circ p$ is a finite type function on $K$. (After all, $f=\Tr(A\pi(\cdot))$ where
$(V,\pi)$ is a  finite dimensional representation of $K/N$ and $A\in\mathrm{End}\,V$. Now
$(V,\pi\circ p)$ is a finite dimensional representation of $K$, and 
$\widehat{f}(k)=\Tr(A(\pi\circ p)(k))$.) As contained in the proof of Lemma \ref{lem-int}, we have
$\int_K \widehat{f}d\mu=\int_{K/N}f\,d\mu'$, where $\mu$ and $\mu'$ are the normalized
Haar measures on $K, K/N$, respectively.

If now $f,g$ are finite type functions on $K/N$ and $\int_{K/N}f^n\,d\mu'=0$ for all $n\in\7N$, then 
$\widehat{f}, \widehat{g}$ are finite type functions on $K$ with
$\int_K \widehat{f}^nd\mu=0\ \forall n$. Since the Mathieu conjecture holds for $K$, we have
\[ \int_{K/N}f^ng\,d\mu'=\int_K\widehat{f^ng}\,d\mu=\int_K \widehat{f}^n \widehat{g}\,d\mu=0\]
for all $n$ large enough. Thus the Mathieu conjecture holds for $K/N$.
\qed

It therefore suffices to prove -- or rather reduce to something simpler -- the Mathieu conjecture for compact
groups of product form. As an immediate corollary of the results of Section \ref{sec5} we have:

\begin{prop}\label{prop-mathieu2}
Let $K, N, M, \delta$ be as in Theorem \ref{thm-main1}. Let
\[ \2A(K)=\{ \stackrel{\approx}{f}\ | \ f\in \7C[K]\}\subset \7C[x_1,\ldots,x_N,z_1^{\pm 1},\ldots,z_M^{\pm 1}].\]
Then the Mathieu conjecture holds for $K$ if and only if for all $f,g\in\2A(K)$ we have that
$\int_{(N,M)} f^n\delta=0$ for all $n\in\7N$ implies
$\int_{(N,M)} f^n g\delta=0$ for all but finitely many $n\in\7N$.
\end{prop}


In this reformulation of the Mathieu conjecture, there is a complicated dependence on the Lie group $K$ both
in the weight function $\delta$ and the structure of the space $\2A(K)$.
Let us again consider the simplest case:

\medskip

\noindent{\bf{Example.}}
For $K=SU(2)$ we have $N=1, M=2$ and $\delta(x)=x$. Taking into account Remark \ref{rem2}.2 and Proposition
\ref{prop-tt} we find that $\2A(SU(2))$ is the set of polynomials in the functions
$w(1-x^2)z,\ xz,\ xz^{-1},\ w^{-1}z^{-1}$, thus
\[ \2A(SU(2)) = \mathrm{span}_\7C\{  x^b (1-x^2)^a w^{a-d} z^{a+b-c-d} \ | \  a,b,c,d\in\7N_0  \},\]
which is a proper subalgebra of $\7C[x,z^{\pm 1}, w^{\pm 1}]$. In the case of $K=SU(2)$, the weight function
$\delta(x)=x$ can be removed by another change of variables, cf.\ \cite{mueger-tus}, but there is no
reason to expect this in general.

\medskip

The main aim of this paper is to formulate conjectures that imply the Mathieu conjecture, all the while being
free of any reference to Lie groups and their complicated structure theory. In view of Propositions
\ref{prop-mathieu1},\ \ref{prop-mathieu2}, the following conjecture clearly implies the Mathieu conjecture: 

\begin{conj}\label{conj-1}
Let $f,g\in\7C[x_1,\ldots,x_N, z_1^{\pm 1},\ldots,z_M^{\pm 1}]$ for non-negative integers $M,N$, and
assume that $\delta\in\7C[x_1,\ldots,x_N]$ is a monomial of odd degree in all variables. If 
$\int_{(N,M)} f^n\delta=0$ for all $n\in\mathbb{N}$ then 
$\int_{(N,M)}f^n g\delta =0$ for all but finitely many $n\in\mathbb{N}$.
\end{conj}

\begin{rema} 1. The special case $N=0$ (in which $\delta\equiv 1$ trivially), equivalent to the Mathieu
conjecture for abelian groups, has been proven in \cite{duist-kallen}. As modest evidence for the general
case, there are so far only the results of \cite{FPYZ} relevant for the case $N=1, M=0$.

2. The pairs $(N,M)$ relevant for the application to the Mathieu conjecture all satisfy $M>N$,
but it seems unlikely that the proof of the conjecture simplifies if $M>N$ is imposed.
\end{rema}

We close by giving a reformulation of Conjecture \ref{conj-1} that was instrumental in \cite{duist-kallen}
(and \cite{mueger-tus}). Namely, with notation and assumptions as in that conjecture, we write
$f=\sum_{{\bf m}\in\mathbb{Z}^M}c_{\bf m} z^{\bf m}$, where $c_{\bf m}\in\7C[x_1,\ldots,x_N]$.
Following \cite{mueger-tus} we define the spectrum of $f$ as
$$
\mathrm{Sp}(f)=\{{\bf m}\in\mathbb{Z}^M\, |\, c_{\bf m}\neq 0\}.
$$

\begin{prop}\label{prop3} 
Keep notation as above. If $\bf 0$ is not in the convex hull of $\mathrm{Sp}(f)$ in $\mathbb{R}^M$, then 
$\int_{(N,M)}f^n g\delta=0$ for all but finitely many $n\in\mathbb{N}$.
\end{prop}

\bp
By a classical result (a baby version of the Hahn-Banach separation theorem), there is a hyperplane in
$\mathbb{R}^M$ separating $\bf 0$ from $\mathrm{Sp}(f)$. This is equivalent to the existence of
${\bf v}\in\mathbb{R}^M\backslash\{{\bf 0}\}$ such that ${\bf v}\cdot{\bf m}\geq 1$ for all
${\bf m}\in\mathrm{Sp}(f)$. Hence ${\bf v}\cdot{\bf m}\geq n$ for all ${\bf m}\in\mathrm{Sp}(f^n)$. So
$\mathrm{Sp}(f^n)$ moves off to infinity as $n\to\infty$, in the sense that it becomes disjoint from every
finite subset of $\mathbb{Z}^M$ for $n$ large enough. In particular, it will eventually be disjoint from
$-\mathrm{Sp}(g)$, showing that the $z$-integration in $\int_{(N,M)}f^ng$ gives zero for all $x\in [0,1]^N$.  
\qed

This motivates a slight generalization of the $xz$-Conjecture from \cite{mueger-tus}:

\begin{conj}\label{conj-2}
Keep notation as above. If $\int_{(N,M)} f^n\delta=0$ for all $n\in\mathbb{N}$, then
$\bf 0$ is not in the convex hull of $\mathrm{Sp}(f)\subset\mathbb{R}^M$.
\end{conj}


Now Proposition \ref{prop3} just says that Conjecture \ref{conj-2} implies Conjecture \ref{conj-1},
thus the Mathieu conjecture.

\begin{rema}
1. It is quite reasonable to expect the following stronger version of Conjecture \ref{conj-2} to hold:
If $\bf 0$ is in the convex hull of $\mathrm{Sp}(f)$ then
$\limsup_{n\rarr\infty}|\int_{(N,M)}f^n\delta|^{1/n}>0$. This statement is true for $N=0$, and its proof in 
\cite{duist-kallen} was the heart of the proof of the abelian Mathieu conjecture given there. (The case
$N=1, M=0, \delta\equiv 1$ was proven in \cite{mueger-tus0}.) 
There is reason to believe that the same will hold for a future proof of the general Conjecture \ref{conj-1}.

2. Since the Mathieu conjecture implies \cite{mathieu} the Jacobian conjecture, the same holds for
the Conjectures \ref{conj-1} and \ref{conj-2}. It would be interesting to (i) find a direct proof of this 
implication, bypassing compact connected Lie groups and Mathieu's conjecture, and to (ii) know whether this
approach to the Jacobian conjecture is related to the known ones.
\end{rema}

\end{document}